\documentstyle{amsart}
\makeatletter
\renewcommand{\subsection}{\@startsection{subsection}{2}{\z@}%
{\baselineskip}{0.5\baselineskip}{\defaultfont\bf}}
\makeatother
\newtheorem{thm}{Theorem}[section]

\newtheorem{lem}[thm]{Lemma}
\theoremstyle{definition}
\newtheorem{defn}{Definition}[section]
\def\dj{d\kern-.30em\raise1.25ex\vbox{\hrule width .3em height .03em}}
\def\Dj{D\kern-.70em\raise0.75ex\vbox{\hrule width .3em height .03em}
\kern.03em}
\newsymbol\rest 1316
\newcommand{\SO}{\mbox{\shape{n}\selectfont SO}}
\newcommand{\U}{\mbox{U}}
\newcommand{\e}{\epsilon}
\newcommand{\k}{\kappa}
\newcommand{\ad}{\mbox{\shape{n}\selectfont ad}}
\newcommand{\id}{\mbox{\shape{n}\selectfont id}}
\newcommand{\vh}{\mbox{\shape{n}\family{euf}\selectfont vh}}
\newcommand{\kr}{\ker}
\newcommand{\grten}{\mathbin{\widehat{\otimes}}}

\newcommand{\lie}{\mbox{\shape{n}\selectfont lie}}
\newcommand{\Ad}{\varpi}
\newcommand{\hor}{\mbox{\family{euf}\shape{n}\selectfont hor}}

\newcommand{\Sum}{\displaystyle{\sum}}
\newcommand{\lc}{\frak{fx}}
\newcommand{\pr}{\pi}
\newcommand{\inv}{i\!\hspace{0.8pt}n\!\hspace{0.6pt}v}
\def\lbl{\star}
\begin{document}
\title{ON FRAMED QUANTUM PRINCIPAL BUNDLES}
\author{Mi\'co \Dj UR\Dj EVI\'c}
\address{Instituto de Matematicas,
UNAM, Area de la Investigacion Cientifica,
Circuito Exterior, Ciudad Universitaria, M\'exico DF, CP 04510,
MEXICO\newline
\indent {\it Written In}\newline
\indent Facultad de Estudios Superiores Cuautitlan, UNAM, MEXICO}
\maketitle
\begin{abstract}
A noncommutative-geometric formalism of framed principal bundles is sketched,
in a special case of quantum bundles (over quantum spaces) possessing
classical structure groups. Quantum counterparts of torsion operators
and Levi-Civita type connections are analyzed. A construction
of a natural differential calculus on framed bundles is described.
Illustrative examples are presented.
\end{abstract}
\section{Introduction}
{\renewcommand{\thepage}{}
     In classical differential geometry the formalism of principal
bundles plays a central role. In particular, a distinguished  role  is
played by {\it framed} principal bundles, characterized as covering
bundles of subbundles of the bundle of linear frames of the base manifold.
Framed  principal  bundles  provide  a   natural   conceptual
framework  for  the  study  of fundamental   classical
differential-geometric structures.

     In this paper classical idea of  a  framed  principal
bundle will be  incorporated  into  an
appropriate  noncommutative-geometric \cite{C} context.
The main entities  figuring  in  the game will  be  quantum
principal bundles (over quantum spaces) possessing (compact)
classical structure groups.

     Conceptually, the paper is based on a general  theory
of quantum principal bundles, presented in \cite{D2}.

     The paper is organized as follows.

     Section II is devoted to the definition and general
properties  of frame structures.
 As  first,  a  very  important  class  of  integrable   frame
structures will be defined.
Roughly  speaking,  integrable  frame  structures   naturally
induce ``Levi-Civita'' type connections on the bundle. On the  other
hand, such connections are compatible, in  the  appropriate  sense,
with internal geometrical structure of the bundle. This fact opens
a possibility to construct,  in  a  fully  intrinsic  manner,  the
complete differential calculus on  the  bundle,  starting  from  a
given integrable frame structure.

     More precisely, starting from a quantum principal bundle
$P$ and an
appropriate representation $u$  of  the  structure  group  $G$  it  is
possible to construct a graded *-algebra $\hor_P$  playing the  role
of horizontal forms on $P$ (together  with  the  right  (co)action  of
$G$ on  this  algebra).  Then,  the  graded  *-algebra
$\Omega_M$
representing differential  forms  on  the  base  manifold  can  be
described   as   the   subalgebra   of $\hor_P$    consisting    of
$G$-invariant elements. On the  other  hand,  if  an  integrable
frame structure on $P$ (with respect to $u$)  is  given,  then  it  is
possible to define, in a natural manner,  a  differential  on  the
algebra
$\Omega_M$. Finally, starting  from  $\hor_P$   and
$\Omega_M$  and
applying ideas of \cite{D3}
it  is  possible  to  construct  (via  the
concept  of  a  preconnection)  the  whole   graded   differential
*-algebra  $\Omega_P$  representing  differential  forms  on  the
bundle $P$.

     In the framework of  (integrable)  frame  structures,  it  is
possible to define  torsion  operators  associated  to  connection
forms (as in the classical theory).
These operators are analyzed at the  end  of Section II.

     It is important to mention that connection forms will  appear
implicitly,   represented   by   the   corresponding    covariant
derivatives. Further, all connections appearing in this study will
be regular and multiplicative (in the sense of \cite{D2})

     In Section III some examples of framed quantum principal bundles
are presented.

     Finally, in Section IV concluding remarks are made.
}
\section{Frame Structures}

     Let $M$ be a quantum space, represented by a *-algebra
${\cal V}$. Let $G$
be an ordinary compact (matrix) Lie group.
Concerning group entities, the notation of \cite{W1}
will be followed (although in the classical context).
In particular the group $G$ will be described by
a (commutative)
Hopf *-algebra $\cal A$ consisting of polynomial functions on $G$.
The group structure is encoded in the coproduct $\phi\colon\cal{A}
\rightarrow\cal{A}\otimes\cal{A}$, counit
$\e\colon\cal{A}\rightarrow\Bbb{C}$ and the antipode $\k\colon\cal{A}
\rightarrow\cal{A}$.

Let $u$ be a real
unitary representation of $G$ in the standard
$n$-dimensional  unitary  space
$\Bbb{C}^n$. We shall assume that the kernel of $u$ is discret.
Furthermore, $u$ will be interpreted as a right comodule structure map
$u\colon\Bbb{C}^n\rightarrow\Bbb{C}^n\otimes\cal{A}$, so that
$$u(e_i)=\sum_{j=1}^ne_j\otimes u_{ji},$$
where $u_{ij}$ are matrix elements of $u$ (and $e_i$ are absolute basis
vectors in $\Bbb{C}^n$).

Let $P=({\cal B},i,F)$ be a  quantum  principal $G$-bundle  over
$M$.
Here, ${\cal  B}$  is a (unital) *-algebra  consisting  of   appropriate
``functions'' on the bundle, $i\colon  {\cal V}\rightarrow{\cal B}$
is the dualized
``projection''  of
$P$ on $M$, and $F\colon  {\cal B}\rightarrow{\cal B}\otimes{\cal  A}$  is
the  dualized
``right action'' of $G$ on $P$. The elements of ${\cal V}$
will be identified with their images from $i({\cal V})$. The algebra
${\cal B}$ is understandable as a  bimodule  over
${\cal  V}$,  in  a natural manner.

 Let us assume that  a  system  $\tau=(\partial_1,\dots,\partial_n)$ of
$\cal B$-valued
hermitian  derivations  $\partial_i\colon  {\cal  V}\rightarrow{\cal
B}$ is given such that
\begin{equation}\label{21}
F\partial_i(f)=\sum_{j=1}^n\partial_j(f)\otimes u_{ji}
\end{equation}
for each $i\in\{1,\dots,n\}$ and $f\in\cal{V}$. Finally,  let  us  assume
that  the  following {\it completeness condition} holds.

\mbox{}

There exist a natural number $d$ and elements
$b_{i\alpha}\in\cal{B}$ and $v_{i\alpha}\in\cal{V}$
(where $\alpha\in\{1,\dots,d\}$ and $i\in\{1,\dots,n\}$)
such that
\begin{equation}\label{22}
\sum_{\alpha}b_{i\alpha}\partial_j(v_{i\alpha})=\delta_{ij}1,
\end{equation}
for each $i,j\in\{1,\dots,n\}$.

\begin{defn} Every system $\tau$ satisfying the
above  conditions
is called a {\it frame structure} on $P$ (relative to $u$).
\end{defn}
\begin{defn} A frame structure $\tau$ is called {\it integrable}
iff there exists  a  system $\widehat{\tau}=(X_1,\dots,X_n)$
of hermitian derivations
$X_i\colon  {\cal B}\rightarrow{\cal B}$ satisfying
\begin{gather}
FX_j=\sum_{k=1}^n(X_k\otimes u_{kj})F\label{23}\\
X_i{\rest}\cal V=\partial_i\label{24}\\
X_i\partial_j-X_j\partial_i=0\label{25}
\end{gather}
for each $i,j\in\{1,\dots,n\}$.
\end{defn}

In the following, it will be assumed  that  the  bundle
$P$ is endowed with a fixed integrable frame  structure
$\tau$. Let us consider a graded *-algebra
$$\hor_P={\cal B}\otimes \Bbb{C}_n^{\wedge},$$
where $\Bbb{C}_n^{\wedge}$  is the
corresponding external *-algebra.
The elements of $\hor_P$ will
be interpreted as ``horizontal forms'' on $P$.
Algebras ${\cal  B}$
and $\Bbb{C}_n^{\wedge}$ are naturally understandable as
subalgebras of $\hor_P$.
We shall denote by
$\theta_i\leftrightarrow 1\otimes e_i$ special horizontal 1-forms
corresponding to absolute basis vectors.
Let $F^\wedge\colon  \hor_P\rightarrow\hor_P
\otimes{\cal A}$ be the product of actions $F$  and $u^\wedge$,
where $u^\wedge\colon \Bbb{C}_n^{\wedge}\rightarrow \Bbb{C}_n^{\wedge}
\otimes{\cal A}$
is the representation of $G$  in  $\Bbb{C}_n^{\wedge}$, induced by $u$.

To each extension $\widehat{\tau}=(X_1,\dots,X_n)$ of $\tau$
it is possible to
associate a first-order antiderivation $\nabla\colon\hor_P
\rightarrow\hor_P$ such that
\begin{gather}
\nabla(b)=\sum_{k=1}^nX_k(b)\theta_k\label{26}\\
\nabla(\theta_i)=0\label{27}
\end{gather}
for each $b\in{\cal B}$ and $i\in\{1,\dots,n\}$. Moreover,
\begin{gather}
\nabla*=*\nabla\label{29}\\
F^\wedge\nabla=(\nabla\otimes\id)F^\wedge\label{210}.
\end{gather}

Let $\lc(P)$ be the set of all antiderivations $\nabla$ constructed
in this way ($\Leftrightarrow$ the set of all ``frame extensions''
$\widehat{\tau}$). Clearly, $\lc(P)$ is a real affine space, in a natural
manner.

Let $\Omega_M\subseteq\hor_P$  be a graded *-subalgebra
consisting of all $F^\wedge$-invariant  elements.
Clearly, $\Omega_M^0={\cal  V}$.
Elements  of  $\Omega_M$ will play the role of
differential forms on $M$. Because of \eqref{210} we have
$$\nabla(\Omega_M)\subseteq\Omega_M,$$
for each $\nabla\in\lc(P)$.

\begin{lem} (i) There exists the common restriction
$d_M\colon  \Omega_M\rightarrow\Omega_M$ of all maps
$\nabla\in\lc(P)$.

(ii) The space $\Omega_M$ is linearly spanned by elements of the form
\begin{equation}\label{lin-gen}
w=f_0d_Mf_1\dots d_Mf_k
\end{equation}
where $f_i\in\cal{V}$.

(iii) We have
\begin{equation}\label{d2=0}
d_M^2=0.
\end{equation}
\end{lem}
\begin{pf}
Let us fix $\nabla\in\lc(P)$ and let $d_M=\nabla{\rest}\Omega_M$.
The completeness condition implies that each element $w\in\hor_P^k$ can
be written in the form
$$ w=\sum_i w_i d_Mf_i$$
where $f_i\in\cal{V}$ and $w_i\in\hor_P^{k-1}$.
This follows from the equality
\begin{equation}\label{Ti}
\theta_i=\sum_{\alpha}b_{i\alpha}d_Mv_{i\alpha}.
\end{equation}
Moreover, without a lack
of generality we can assume that $w_i\in\Omega_M^{k-1}$. Therefore the
statement follows by applying the principle of mathematical induction.

It is sufficient to check that
\eqref{d2=0} holds on elements of the form \eqref{lin-gen}. Because of
the graded Leibniz rule it is sufficient to check that
$d_M^2(\cal{V})=\{0\}$. However, this directly follows from \eqref{25}
and from the definition of $\nabla$.

Finally, ({\it i}\/) follows from the graded Leibniz rule,
({\it ii}\/) and ({\it iii}\/),
and from the fact that maps from  $\lc(P)$ act on
elements from $\cal{V}$ in the same way (fixed by $\tau$).
\end{pf}

In other words $\Omega_M$, endowed with $d_M$,
becomes a graded-differential *-algebra generated by $\cal{V}$.

     All basic structural elements of the conceptual framework  of
\cite{D3}  are  now  in  the game.
Let us recall that a  {\it  preconnection}  on  $P$  (relative  to
$\bigl\{\hor_P,F^\wedge,\Omega_M\bigr\}$)
is a first-order hermitian
antiderivation $D\colon\hor_P\rightarrow\hor_P$ satisfying
\begin{gather}
F^\wedge D=(D\otimes \id)F^\wedge\label{218}\\
D{\rest}\Omega_M=d_M\label{220}
\end{gather}
(according to \eqref{218} every $D$ is reduced in $\Omega_M$).
Preconnections form  a  real  affine  space  $\pr(P)$.  Evidently,
$\lc(P)\subseteq\pr(P)$.

According to \cite{D3}, for each $D\in\pr(P)$ and $E\in
\overrightarrow{\pr}(P)$
(the vector space associated to $\pr(P)$) there exists the unique
linear maps $\varrho_D^\lbl, \chi_E^\lbl\colon\cal{A}\rightarrow\hor_P$ such
that
\begin{gather}
D^2(\varphi)=-\sum_k\varphi_k\varrho_D^\lbl(c_k),\label{222}\\
E(\varphi)=-(-)^{\partial\varphi}\sum_k\varphi_k \chi_E^\lbl(c_k)
\label{226}
\end{gather}
for each $\varphi\in\hor_P$, where
$\Sum_k  \varphi_k  \otimes  c_k  =F^\wedge(\varphi)$.
The maps $\varrho_D^\lbl$ and $\chi_E^\lbl$ determine,
in a natural manner, a bicovariant first-order *-calculus $\Psi$ on $G$.
This calculus is based (in the sense of \cite{W2}) on the right
$\cal{A}$-ideal
$\cal{R}\subseteq\kr(\e)$ consisting of elements anihilated by
all $\varrho_D^\lbl$ and $\chi_E^\lbl$. These maps can
be therefore factorized through
$\cal{R}$. In such a way we obtain maps $\varrho_D, \chi_E\colon\Psi_{\inv}
\rightarrow\hor_P$, where $\Psi_{\inv}=\kr(\e)/\cal{R}$ is the space of
left-invariant elements of $\Psi$. More precisely,
\begin{gather*}
\varrho_D\pi=\varrho_D^\lbl\\
\chi_E\pi=\chi_E^\lbl
\end{gather*}
where
$\pi\colon\cal{A}\rightarrow\Psi_{\inv}$ is the canonical projection map.

The following equalities hold
\begin{align}
F^\wedge\varrho_D&=(\varrho_D\otimes \id)\Ad\label{223}\\
D\varrho_D&=0\label{224}\\
F^\wedge\chi_E&=(\chi_E\otimes \id)\Ad\label{227}
\end{align}
where $\Ad\colon \Psi_{\inv}\rightarrow \Psi_{\inv}\otimes\cal{A}$ is the
dualized (co)adjoint action, explicitly given by
$$\Ad\pi=(\pi\otimes\id)\ad,$$
and $\ad\colon\cal{A}\rightarrow\cal{A}\otimes\cal{A}$ is the dualized
adjoint action of $G$ on itself. Further,
\begin{align}
\chi_E(\vartheta)\varphi&=(-)^{\partial\varphi}
\sum_k\varphi_k\chi_E(\vartheta\circ c_k)\label{228}\\
\varrho_D(\vartheta)\varphi&=\sum_k\varphi_k\varrho_D
(\vartheta\circ c_k)\label{225}
\end{align}
for each $\vartheta\in \Psi_{\inv}$
and $\varphi\in\hor_P$, where $\circ$ is the canonical right
$\cal{A}$-module structure.

The map $\varrho_D$  is called {\it the curvature} of $D$.

Every  element
$\varphi\in\hor_P$ can be written in the form
$$\varphi=\sum_kb_kd_Mw_k, $$
where $b_k \in{\cal B}$ and $w_k\in\Omega_M$.
This implies that each
$D\in\pr(P)$ is completely determined  by  its  restriction
$D{\rest}{\cal  B}$.
Explicitly, this restriction is described by
\begin{equation}\label{D=SY}
D(b)=\sum_{i=1}^n Y_i(b)\theta_i,\label{229}
\end{equation}
where $Y_i\colon\cal{B}\rightarrow\cal{B}$ are hermitian
derivations satisfying
\begin{align}
Y_i{\rest}\cal V&=\partial_i\label{230}\\
FY_i(b)&=\sum_{kj}Y_j(b_k)\otimes c_k u_{ji}\label{231}
\end{align}
for each $b\in{\cal B}$, where $\Sum_k b_k \otimes c_k=F(b)$. In
terms of derivations $Y_i$ the action of $D$ on special horizontal forms
$\theta_j$ is given by
\begin{equation}\label{DT}
D(\theta_j)=\frac{1}{2}\sum_{\alpha kl}\Bigl\{
Y_k(b_{j\alpha})\partial_l(v_{j\alpha})-Y_l(b_{j\alpha})\partial_k(v_{j
\alpha})\Bigr\}\theta_k\theta_l.
\end{equation}

For every $D\in \pr(P)$ let $\Theta_D\colon \Bbb{C}^n\rightarrow\hor_P$
be a linear map given by
\begin{equation}
\Theta_D(e_i)=\Theta_D^i=D(\theta_i).\label{232}
\end{equation}
\begin{defn}
The map $\Theta_D$ is called {\it the torsion} of $D$.
\end{defn}
\begin{lem} The map $\Theta_D$ is hermitian and
satisfies
\begin{align}
F^\wedge\Theta^i_D&=\sum_{j=1}^n\Theta_D^j\otimes u_{ji}\label{233}\\
-D\Theta_D^i&=\sum_{j=1}^n \theta_j \varrho_D^\lbl(u_{ji})\label{234}
\end{align}
for each $i\in\{1,\dots,n\}.$
\end{lem}
\begin{pf}
The statement follows directly from  properties
\eqref{218} and \eqref{222} and from the
definition of $\Theta_D$ and $F^\wedge$.
\end{pf}

     Identity \eqref{234} corresponds to the second Structure  equation
in  classical  differential  geometry  \cite{KN}.

\begin{lem}
The following equivalence holds
\begin{equation}
\Theta_D=0\iff D\in\lc(P)
\end{equation}
for each $D\in\pr(P)$.
\end{lem}
\begin{pf}
Let us assume that $D\in\lc(P)$. From \eqref{22} and
\eqref{DT} it follows that the torsion vanishes. Conversely,
if $\Theta_D=0$ then
$$0=d_M^2(f)=D^2(f)=\frac{1}{2}\sum_{ij}[Y_i,Y_j](f)\theta_i\theta_j$$
for each $f\in\cal{V}$. In other words, $D\in\lc(P)$.
\end{pf}

Let us assume that the higher-order differential calculus on $G$ is described
by the universal envelope $\Psi^\wedge$ of $\Psi$.
Let $\Omega_P$ be the graded-differential *-algebra
canonically associated, in the sense of \cite{D3}, to algebras
$\hor_P$  and
$\Psi^\wedge$, and to the system $\pr(P)$ of preconnections. The main
property of this algebra is that each $D\in\pr(P)$ naturally
induces a representation of the form
\begin{equation}\label{dec-W}
\Omega_P\leftrightarrow\vh_P
\end{equation}
where $\vh_P$ is a graded *-algebra representing
``vertically-horizontally'' decomposed forms \cite{D2} on $P$. At
the level of graded
vector spaces, we have
$$\vh_P=\hor_P\otimes\Psi_{\inv}^\wedge$$
The *-algebra structure on
$\vh_P$ is specified by
\begin{gather*}
(\psi\otimes\vartheta)(\varphi\otimes\eta)=\sum_k(-)^{\partial
\vartheta\partial\varphi}\psi\varphi_k\otimes(\vartheta\circ c_k)\eta\\
(\varphi\otimes\vartheta)^*=\sum_k\varphi_k^*
\otimes(\vartheta^*\circ c_k^*).
\end{gather*}
In terms of the identification \eqref{dec-W} the differential structure
on $\Omega_P$ is expressed
via a $D$-dependent differential $\partial_D$ on $\vh(P)$, which is given by
\begin{gather*}
\partial_D(\varphi)=D(\varphi)+(-)^{\partial\varphi}\sum_k\varphi_k\pi(c_k)\\
\partial_D(\vartheta)=\varrho_D(\vartheta)+d(\vartheta)
\end{gather*}
where $\vartheta\in\Psi_{\inv}$ and $d\colon\Psi_{\inv}^\wedge
\rightarrow\Psi_{\inv}^\wedge$ is the corresponding differential
($\partial_D$ is extended on $\vh_P$ by the graded Leibniz rule).

\begin{lem}
As  a differential   algebra, $\Omega_P$ is   generated   by
${\cal B}=\Omega^0_P$.
\end{lem}
\begin{pf}
Because of \eqref{Ti}, it is sufficient to check that elements of the
form $\pi(a)$ belong to $\cal{B}\partial_D\bigl\{\cal{B}\bigr\}$.
For a given $a\in\cal{A}$, let us choose elements
$q_k,b_k\in\cal{B}$ such that $\Sum_k q_kF(b_k)=1\otimes a$ (the group $G$
acts freely on $P$).
Then the following equality holds
$$\sum_k q_k\partial_D(b_k)=\sum_kq_kD(b_k)+\pi(a).$$
This implies (together with \eqref{Ti} and \eqref{D=SY})
that $\pi(a)$ is expressible in the desired way.
\end{pf}

The  map  $F$  is   uniquely   extendible
to the homomorphism
$\widehat{F}\colon  \Omega_P\rightarrow\Omega_P\grten
\Gamma^\wedge$ of graded differential *-algebras (corresponding
to the ``pull back'' map of differential forms).

     There exists  a  natural  bijective  (affine)  correspondence
$D\leftrightarrow\omega$ between preconnections  $D$  and  regular
connections $\omega$
on $P$. In terms of this correspondence,
\begin{gather*}
R_\omega \leftrightarrow\varrho_D\\
D_\omega\leftrightarrow D.
\end{gather*}

It is also possible to construct differential structures
on $G$ and $P$ starting from a restricted set of preconnections
forming an affine subspace of $\lc(P)$. In this case covariant derivatives
of regular connections will (generally) form only an affine subspace of
$\pr(P)$.
In particular,
every {\it single} element
$\nabla\in\lc(P)$ determines a differential calculus on $P$. Let
us assume that $\Psi$ is the minimal bicovariant first-order
differential *-calculus over $G$ compatible with $\nabla$. Then the
corresponding caclulus on $P$ will be based on
a graded-differential *-algebra
$\Omega_P=(\vh_P,\partial_\nabla)$.
It is also possible to vary this theme,
and to choose for $\Psi$ an arbitrary (non-minimal)
calculus satisfying
the mentioned compatibility conditions.

Let $\pr(P)_\nabla\subseteq\pr(P)$ be the affine subspace consisting
of preconnections interpretable as covariant derivatives of
regular connections (relative to $\Omega_P$ constructed from
$\bigl\{\nabla\bigr\}$).

A particularly interesting situation arises when
$\nabla$ is compatible
with
{\it the classical} differential calculus on $G$.
\begin{defn}
An element $\nabla\in\lc(P)$ is called {\it classical} iff
$$\varrho_\nabla^\lbl(ab)=\epsilon(a)\varrho_\nabla^\lbl(b)+
\varrho_\nabla^\lbl(a)\epsilon(b)$$
for each $a,b\in\cal{A}$.
\end{defn}

If $\nabla$ is classical then it is possible to
assume that $\Psi$ is the classical differential calculus
on $G$ (hence $\Psi_{\inv}=\lie(G)^*$).
In this particular case,
\begin{align*}
\varrho_D^\lbl(ab)&=\epsilon(a)\varrho_D^\lbl(b)+
\varrho_D^\lbl(a)\epsilon(b)\\
\chi_E^\lbl(ab)&=\epsilon(a)\chi_E^\lbl(b)+
\chi_E^\lbl(a)\epsilon(b)\\
\varrho_D(\vartheta)\varphi&=\varphi\varrho_D(\vartheta)\\
\chi_E(\vartheta)\varphi&
=(-)^{\partial\varphi}\varphi\chi_D(\vartheta)
\end{align*}
for each $D\in\pr(P)_\nabla$ and $E\in\overrightarrow\pr(P)_\nabla$
(the $\circ$-structure on $\Psi_{\inv}$ is trivialized, because
$\vartheta\circ a=\epsilon(a)\vartheta$).

\begin{lem}We have
\begin{equation}
\pr(P)_\nabla\cap\lc(P)=\bigl\{\nabla\bigr\}
\end{equation}
for each $\nabla\in\lc(P)$.
\end{lem}
\begin{pf}
Let us consider elements $D\in\pr(P)_\nabla$ and $\nabla\in\lc(P)$, and let
us assume that
$$
\nabla\leftrightarrow (X_1,\dots,X_n)\qquad
D\leftrightarrow (Y_1,\dots,Y_n).
$$
Derivations $Z_i=X_i-Y_i$ possess the following properties
\begin{align}
FZ_i(b)&=\sum_{kj}Z_j(b_k)\otimes c_k u_{ji}\label{235}\\
Z_i(f)&=0\label{236}
\end{align}
for each $b\in{\cal B}$ and $f\in{\cal V}.$
Applying results of \cite{D3}
we conclude that there exist linear maps
$\lambda_i\colon \Psi_{\inv}\rightarrow{\cal B}$ such that
\begin{equation}
Z_i(b)=\sum_k b_k\lambda_i\pi(c_k)\label{238}
\end{equation}
for each $b\in{\cal B}$ and $i\in\{1,\dots,n\}.$
In particular,
\begin{equation}
Z_i\partial_j(f)=\sum_k \partial_k (f)\lambda^k_{ij}\label{239}
\end{equation}
where
\begin{equation}
\lambda_{ij}^k=\lambda_i\bigl(\pi(u_{kj})\bigr).\label{240}
\end{equation}
Let us assume that $\Theta_D=0$. This is equivalent to
\begin{equation}
(Z_i\partial_j -Z_j\partial_i)(f)=0\label{237}
\end{equation}
for each $f\in\cal V$ and $i,j\in\{1,\dots,n\}$.
Identities \eqref{239}--\eqref{237}, together with the  completeness
condition imply
\begin{equation}
\lambda_{ij}^k=\lambda_{ji}^k\label{241}
\end{equation}
On the other hand
\begin{equation}
\lambda_{ij}^k=-\lambda_{ik}^j\label{242}
\end{equation}
as easily follows from \eqref{240} and the hermicity of $u_{ij}$.
It follows that $\lambda_{ij}^k=0$, and hence $\lambda_i=0$,
because $\Psi_{\inv}$ is
spanned by elements $\pi(u_{ij})$. Hence $D=\nabla$.
\end{pf}

The above equivalence
corresponds  to  the  classical  characterization   of the
Levi-Civita connection, as the unique  (metric)  connection  with
vanishing torsion. Conceptually, we  followed  a  classical  proof
\cite{KN} of the uniqueness of the Levi-Civita connection.

\section{Examples}
\subsection{The classical case}
Let $P$ be a classical principal $\SO(k)$-bundle  over
a compact smooth $n$-dimensional maniflold $M$  (where
$k\leq n$) and let $\tau=(\partial_1,\dots,\partial_k)$
be a frame structure
on $P$ (relative to the standard representation $u$ of $\SO(k)$ in
$\Bbb{C}^k$).
Then every point $p\in P$ naturally determines a $k$-tuple
$(\xi_1,\dots,\xi_k)$
on tangent vectors on $M$ in the point $x=\pi_M(p)$ as follows
\begin{equation}
\xi_i(f)=[\partial_i(f)](p).\label{32}
\end{equation}
Here $\pi_M\colon   P\rightarrow M$ is the projection map.

     From the transformation property \eqref{21},
it follows  that  the
space $\Sigma_x\subseteq T_x(M)$ spanned by $(\xi_1,\dots,\xi_k)$
is independent
of the choice of the point $p\in \pi_M^{-1}(x)$.
On  the  other  hand, completeness
condition \eqref{22} implies that $(\xi_1,\dots,\xi_k)$  are
linearly independent vectors.

     In such a way an oriented $k$-dimensional subbundle $\Sigma$  of
$T(M)$
is constructed. Fibers  $(\Sigma_x)_{x\in M}$ possess a natural  Euclidean
structure defined by requiring that  $(\xi_1,\dots,\xi_k)$
are orthonormal vectors.
     In classical terms, $P$ is identificable  with  the  bundle  of
oriented orthonormal frames of $\Sigma$.

Let us assume that $\tau$ is integrable. This implies that the space of
smooth sections
of $\Sigma$ is closed with respect to the commutator of vector fields.
In other words, $\Sigma$ is integrable (according to Frobenius  theorem).

     Let $N\subseteq M$ be an arbitrary leaf of the foliation
$\Sigma$, and
$P_N$   the
portion of $P$ over $N$. This bundle  coincides, in a natural manner,
with
the  bundle  of
oriented orthonormal frames of $N$. It is invariant under the action
of fields $X_i$. Restrictions of $X_i$  on $P_N$ determine
the standard Levi-Civita connection.

     Bundles $P_N$   determine  a  foliation  $\Sigma^\lbl$ of $P$.
The  elements of $\Omega_P$ are naturally interpretable as
$\Sigma^\lbl$-differential forms  on $P$. In this picture
the algebra $\Omega_M$ consists of $\Sigma$-differential
forms on $M$.

In particular, the case $k=n$ is equivalent to
the classical oriented Riemannian manifold structure on $M$, so that
$P$ becomes the corresponding bundle of oriented orthonormal frames.

\subsection{A Framed Quantum $\SO(2)$ Bundle}

Let us assume that $\cal{V}$ is endowed with a *-automorphism
$\gamma\colon\cal{V}\rightarrow\cal{V}$. Let us assume that
$G=\SO(2)$. The Hopf *-algebra $\cal{A}$ of polynomial functions
on $G$ is generated by a unitary element $U=\cos+i\sin$ (and $\{
\cos, \sin\}$ are understood as functions on $G$). We have
$\phi(U)=U\otimes U$. The formulas
\begin{align}
(f\otimes U^m)(g\otimes U^n)&=f\gamma^m(g)\otimes U^{m+n}\\
(f\otimes U^m)^*&=\gamma^{-m}(F^\wedge)\otimes U^{-m}
\end{align}
(where $n,m\in\Bbb{Z}$)
define a *-algebra structure on the vector space $\cal{B}=
\cal{V}\otimes\cal{A}$. Let $i\colon\cal{V}\rightarrow
\cal{B}$ be the canonical inclusion map. The formula
$$ F(f\otimes U^m)=f\otimes U^m\otimes U^m$$
defines an action of $G$ by *-automorphisms of $\cal{B}$, so that
$P=(\cal{B},i,F)$ is a (quantum) principal $G$-bundle over $M$.

Let us define derivations $X_\pm\colon\cal{B}\rightarrow
\cal{B}$ by
\begin{gather}
X_+(b)=(\alpha\otimes U)b-b(\alpha\otimes U)\\
X_-(b)=(\beta\otimes \bar{U})b-b(\beta\otimes \bar{U})
\end{gather}
where $\alpha,\beta\in\cal{V}$ are such that
$\beta=-\gamma^{-1}(\alpha^*)$. We have
$$[X_+,X_-](b)=vb-bv$$
where $v=\alpha\gamma(\beta)-\beta\gamma^{-1}(\alpha)$. Let
$X_1,X_2\colon\cal{B}\rightarrow\cal{B}$ be (hermitian) derivations
given by $X_\pm=X_1\mp iX_2$. Let $\partial_1,\partial_2=
X_1,X_2{\rest}\cal{V}$ be the corresponding
restrictions. If the completeness condition \eqref{22}
holds then
$\tau=(\partial_1,\partial_2)$ is a frame structure on $P$, relative
to the standard representation
\begin{equation*}u=
\begin{pmatrix}
\cos&-\sin\\
\sin&\phantom{-}\cos
\end{pmatrix}
\end{equation*}
of $\SO(2)$ in $\Bbb{C}^2$.
If $v$ is a central element of $\cal{V}$ then $\tau$ is integrable,
and we can write $\widehat{\tau}=(X_1,X_2)$. The corresponding curvature
is given by
\begin{equation}
\varrho^\lbl_\nabla(U^m)=\frac{1}{4i}\bigl(v-\gamma^{-m}(v)\bigr)
\theta_1\theta_2.
\end{equation}
In general, such frame structures induce nonstandard differential calculi
on $G$. Let us assume that $v$ is non-trivial, and that $\gamma(v)=tv$, for
some $t\in\Re\setminus\{-1,0,1\}$.
This naturally induces a $1$-dimensional calculus on $G$.

The space
$\Psi_{\inv}$ is spanned by $\zeta=\pi(U-\bar{U})$. The corresponding ideal
$\cal{R}$ is generated by $tU+\bar{U}-(1+t)1$. The $\circ$-structure
on $\Psi_{\inv}$
is specified by $\zeta\circ U^m=t^{-m}\zeta$.

\subsection{Free Actions of Simple Lie Groups on Quantum Spaces}
Let us assume that  a  compact  simple  Lie  group
$H$ acts freely on a quantum  space $P$  (determined  by  a
*-algebra $\cal{B}$). Let $G$  be  a  compact  subgroup   of $H$.
Let
$F^\lbl\colon  {\cal B}\rightarrow{\cal B}\otimes{\cal  A}^\lbl$
be the dualized action of $H$ on $P$
(where  ${\cal  A}^\lbl$  is  the  *-algebra  of
polynomial
functions on $H$) and    let
$F=(\id\otimes j)F^\lbl$
be the restriction of the action of $H$ on $G$. Here $j\colon\cal{A}^\lbl
\rightarrow\cal{A}$ is the restriction map. In what follows
the entities endowed with $\lbl$ will refer to $H$.
Both groups will be endowed with standard differential structures.

The triplet $P=({\cal B},i,F)$
is a quantum principal $G$-bundle over $M$, where $M$  is  the  quantum
space based on  the  $F$-fixed-point  *-subalgebra  ${\cal  V}$,  and
$i\colon  {\cal V}\hookrightarrow{\cal B}$ is the inclusion map.

Let $\frak{g}^\lbl$ be the (complex) Lie
algebra of $H$.
Let $u$ be the adjoint representation of $G$ in the space
$\frak{g}^\perp\subseteq\frak{g}^\lbl$.
Here, $\frak{g}\subseteq \frak{g}^\lbl$
is the Lie algebra of $G$, and it is  assumed  that
$\frak{g}^\lbl$ is
endowed with the (positive) Killing scalar product. Let us  assume
that the kernel of $u$ is discrete. Furter, let us assume that
$\frak{g}^\perp$ is identified with $\Bbb{C}^n$, with the help of
a real orthonormal basis $(\xi_1,\dots,\xi_n)$ in $\frak{g}^\perp$.

For each $i\in\{1,\dots,n\}$ let $X_i=(\id\otimes\xi_i\pi^\lbl)F^\lbl$  be
the  hermitian
derivation on ${\cal B}$ corresponding to $\xi_i$ (here we have
identified ${\vphantom{\Psi}\frak{g}}^\lbl=(\Psi_{\inv}^\lbl)^*$)

Let $\partial_i\colon  {\cal V}\rightarrow{\cal B}$
be the restrictions of $X_i$ on ${\cal V}$.
\begin{lem}
Under      the      above      assumptions
$\tau=(\partial_1,\dots,\partial_n)$ is a frame  structure  on
$P$.  This structure is integrable if
\begin{equation}
[\frak{g}^\perp,\frak{g}^\perp]\subseteq\frak{g}.\label{312}
\end{equation}
In this case
\begin{equation}\label{313}
\widehat{\tau}=(X_{1},\dots,X_{n}).
\end{equation}
\end{lem}
\begin{pf}
A direct computation gives
\begin{equation*}
\begin{split}
FX_{i}(b)&=\sum_k F\bigl(b_k \xi_i\pi^\lbl(d_k)\bigr)\\
&=\sum_k b_k \xi_i
\pi^\lbl(d_k^{(2)})\otimes j\left[d_k^{(1)}\k^\lbl(d_k^{(3)})
d_k^{(4)}\right]
=\sum_{jk}X_j(b_k)\otimes u_{ji}j(d_k).
\end{split}
\end{equation*}
Here, $F^\lbl(b)=\Sum_k b_k\otimes d_k$  and we
have applied the identity
$$\xi_i\pi^\lbl(a^{(2)})\otimes j\bigl(a^{(1)}\k^\lbl(a^{(3)})\bigr)
=\sum_j \xi_j
\pi^\lbl(a)\otimes u_{ji}.$$
Hence, derivations $X_i$ transform in the  appropriate  way.

Let  us
prove the completeness condition.
It is sufficient to prove that
for each $a\in\cal{A}^\lbl$, invariant under the right action of
$G$ on $H$, there
exist elements $b_k\in{\cal B}$ and
$v_k \in{\cal V}$ such that
\begin{equation}
\sum_k b_k F^\lbl(v_k)=1\otimes a.\label{315}
\end{equation}
Indeed,  for  a  given $G$-invariant element $a\in\cal{A}^\lbl$
there   exist   elements
$q_k,b_k\in{\cal B}$ such that
$$\sum_k b_k F^\lbl(q_k)=1\otimes a$$
(this is the place where the freeness assumption enters the game).
This imples
$$ \sum_{kl}b_k q_{kl}\otimes
a_{kl}^{(1)}hj(a_{kl}^{(2)})=1\otimes a$$
where  $F^\lbl(q_k)=\Sum_l q_{kl}  \otimes a_{kl}$ and
$h\colon  {\cal  A}\rightarrow \Bbb{C}$
is  the  Haar
measure on $G$. If we define
$$v_k=\sum_l q_{kl}hj(a_{kl}),$$
then $v_k\in{\cal V}$ and \eqref{315} holds.

     Finally, let us assume that \eqref{312} holds. Then
$[X_i,X_j](v)=0$
for each $v\in\cal{V}$ and $i,j\in\{1,\dots,n\}$,
according to the definition of
${\cal V}$. Hence, $\tau$ is integrable and \eqref{313} holds.
\end{pf}

     The corresponding graded-differential *-algebra
$\Omega_P$ can be naturally realized as
$$\Omega_P=[\Psi_{\inv}^\lbl]^\wedge\otimes\cal{B},$$
in other words,
${\cal  B}$-valued  forms  on  $\frak{g}^\lbl$ (with  the
standard   algebraic
structure).

The curvature map of the Levi-Civita connection $\nabla$ is
given by
$$ \varrho_\nabla^\lbl(a)=-\frac{1}{2}\sum_{ij}[\xi_i,\xi_j]\pi(a)
\theta_i\theta_j.$$

\section{Concluding remarks}
     Quantum  counterparts  of  various  important
differential-geometric structures can be introduced in
the  framework  of  the
concept of the frame structure. In particular, frame structures on
quantum $\SO(n)$-bundles provide  a  natural  framework  for  a
noncommutative-geometric version of (oriented)  Riemannian  geometry
(the Xodge *-operator and the  Laplace  operator,  for  example).

     All essential elements of the algebraic  structure  appearing
in the theory of Kahler manifolds are preserved  in  its
noncommutative-geometric version, dealing with framed quantum
$\U(n)$-bundles. In particular, quantum spaces $M$  considered  in  the
previous section can be naturally  endowed  with  Kahler  manifold
structures, in accordance with the analogy
$M \leftrightarrow \mbox{CP}(n)$.

 In this paper, we have assumed that the structure  quantum  group
is compact. This assumption is not essential.  The  whole
formalism  can  be  directly  incorporated  in  a   more   general
conceptual framework, including non-compact structure groups
and  non-unitary representations $u$.
However in this case
(as in classical geometry) $\nabla$ is generally
not uniquely determined by its class $\pr(P)_\nabla$,
even if the calculus on the structure
group is classical.

In terms of the constructed differential calculus on $P$, frame
structures are completely represented by $n$-tuples
$\theta=(\theta_1,\dots,\theta_n)$ of special horizontal 1-forms
$\theta_i$ (counterparts
of classical frame forms \cite{KN}). These forms transform
covariantly, according to the representation $u$.

In this sense,
integrable frame structures can be viewed as
a special case of frame structures introduced in \cite{PH}--Appendix B.

The presented formalism admits a  natural  generalization  to
the fully quantum context (in which the structure group $G$  is
a quantum object). The only essentially new phenomena naturally appearing
in the game is
the presence of a non-trivial right $\cal{A}$-module structure
$\circ$ on the representation space $\Bbb{C}^n$, which is compatible with
the right comodule structure $u$ in such a way
that $\Bbb{C}^n$ together with $u$ and $\circ$ becomes a
``left-invariant part'' of a bicovariant \cite{W2}
*-bimodule $\Gamma$ over $G$. This requires the following compatibility
condition
$$u(e_i\circ a)=\sum_j (e_j\circ a^{(2)})\otimes\k(a^{(1)})u_{ji}a^{(3)}$$
for each $i\in\{1,\dots,n\}$ and $a\in\cal{A}$.

The external algebra $\Bbb{C}_n^\wedge$ should be replaced by
the quantum external algebra, associated to the ``left-invariant part''
$\sigma
\colon\Bbb{C}^n\otimes\Bbb{C}^n\rightarrow\Bbb{C}^n\otimes\Bbb{C}^n$
of
the corresponding canonical flip-over \cite{W2} operator (or its scalar
multiple).
Also, the construction of the horizontal algebra incorporates
elements of the bicovariant bimodule structure.
Maps $X_i$ are not derivations, although their restrictions
$\partial_i=X_i{\rest}\cal{V}$
still satisfy the Leibniz rule.
\filbreak

\end{document}